# A Computer Verification of a Conjecture About The Erdös-Mordell Curve

*Bojan D. Banjac, Branko J. Malešević, Maja M. Petrović and Marija Đ. Obradović*

*Abstract* — *In this paper we consider Erdös-Mordell inequality and its extension in the plane of triangle to the Erdös-Mordell curve. Algebraic equation of this curve is derived, and using modern computer tools in mathematics, we verified one conjecture that relates to Erdös-Mordell curve.*

## I. INTRODUCTION

IN theory of geometrical inequalities [1]-[4], a special place holds the Erdös-Mordell inequality for triangle:

$$R_A + R_B + R_C \geq 2(r_a + r_b + r_c) \qquad (1)$$

where $R_A$, $R_B$ and $R_C$ are the distances from the arbitrary point $M$ in the interior of $\Delta ABC$ to the vertices $A$, $B$ and $C$ respectively, and $r_a$, $r_b$ and $r_c$ are the distances from the point $M$ to the sides $BC$, $CA$ and $AB$ respectively (Fig. 1).

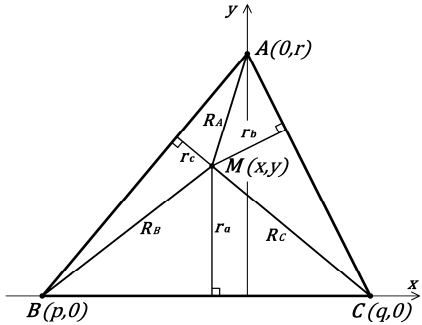

Fig. 1. Erdös-Mordell inequality.

The importance of the Erdös-Mordell inequality can be observed in a recent result of V. Pambuccian which proved that this inequality is an equivalent of negative curvature of absolute plane [5]. In the paper [6], is given an extension of the Erdös-Mordell inequality from the interior of triangle to the interior of the Erdös-Mordell curve. Also, the Erdös-Mordell inequality has its interpretation in the following geometric problem: "Which set of paths/communication lines is the fastest/shortest?" [7].

In this paper we give the algebraic equation of the Erdös-Mordell curve, and in addition, in the case of equilateral triangle, we give a conjecture about the Erdös-Mordell curve brought out from [6], which is verified on the large number of triangles. Value of constant $\varepsilon_0$ from the conjecture is determined by use of symbolic-numerical functions of MatLAB and the conjecture is numerically tested in programming language Java with the use of visual representation of equidistant region in $\varepsilon \geq \varepsilon_0$.

## II. ERDÖS-MORDELL CURVE

In paper [6], a set of points E was introduced, which fulfils:

$$R_A + R_B + R_C \geq \left(\frac{c}{b}+\frac{b}{c}\right)r_a + \left(\frac{c}{a}+\frac{a}{c}\right)r_b + \left(\frac{a}{b}+\frac{b}{a}\right)r_c \qquad (2)$$

where $a=|BC|$, $b=|CA|$ and $c=|AB|$. The set $E$ is determined in the intersection of areas:

$$E_A = \left\{(x,y) \,\Big|\, R_A \geq \frac{c}{a}r_b + \frac{b}{a}r_c\right\}, \qquad (3)$$

$$E_B = \left\{(x,y) \,\Big|\, R_B \geq \frac{c}{b}r_a + \frac{a}{b}r_c\right\}, \qquad (4)$$

$$E_C = \left\{(x,y) \,\Big|\, R_C \geq \frac{b}{c}r_a + \frac{a}{c}r_b\right\}. \qquad (5)$$

The previous sets which correspond to the points $A$, $B$ and $C$ are corner areas, see [6]. If (2) is true, then the Erdös-Mordell inequality is also fulfilled. In the paper [6], the Erdös-Mordell curve is defined by the following equation:

$$R_A + R_B + R_C = 2(r_a + r_b + r_c), \qquad (6)$$

where

$$R_A = \sqrt{x^2 + (y-r)^2}, \qquad r_a = \frac{|y(q-p)|}{\sqrt{(q-p)^2}} = |y|,$$

$$R_B = \sqrt{(x-p)^2 + y^2}, \qquad r_b = \frac{|-q(y-r)-rx|}{\sqrt{r^2+q^2}},$$

$$R_C = \sqrt{(x-q)^2 + y^2}, \qquad r_c = \frac{|-p(y-r)-rx|}{\sqrt{r^2+p^2}},$$

where $A=A(0,r)$, $B=B(p,0)$ and $C=C(q,0)$. Let us denote $E'$ as interior of the Erdös-Mordell curve; then the set $E'$ is maximal set of points in the plane where the Erdös-Mordell inequality is true. In the paper [6] it is proven that the set $E'$ contains the set $E$, and the set $E$ contains initial triangle $\Delta ABC$. The curve (6) is a union of parts of algebraic curves of eight order (Fig. 2), which we are proving in this section.

Bojan D. Banjac is with the School of Electrical Engineering, University of Belgrade, Bulevar kralja Aleksandra 73, 11120 Belgrade, Serbia (phone: 381-64-6356346; e-mail: bojan.b.mail@gmail.com).

Corresponding author Branko J. Malešević, is now with the School of Electrical Engineering, University of Belgrade, Bulevar kralja Aleksandra 73, 11120 Belgrade, Serbia (phone: 381-11-3218317; e-mail: malesevic@etf.rs).

Maja M. Petrović is with the Faculty of Transport and Traffic Engineering, University of Belgrade,Vojvode Stepe 305, 11000 Belgrade, Serbia (phone:381-11-3091259; e-mail: majapet@sf.bg.ac.rs).

Marija Đ. Obradović is with the Faculty of Civil Engineering, University of Belgrade, Bulevar kralja Aleksandra 73 11000 Belgrade, Serbia (phone: 381-11-3218752; e mail: marijao@grf.bg.ac.rs).

Acknowledgment. Research is partially supported by the Ministry of Science and Education of the Republic of Serbia, Grant No. III 44006 and ON 174032.

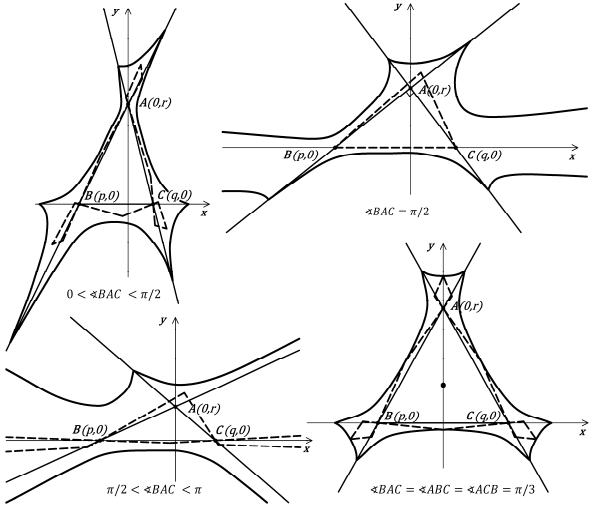

Fig. 2. Erdös-Mordell curve and the area E.

We start from the Erdös-Mordell curve defined by equality (6). Let us denote expressions:
$Q_1 = x^2 + (y-r)^2, Q_2 = (x-p)^2 + y^2, Q_3 = (x-q)^2 + y^2$
and denote a sum of absolute values
$$S = 2(r_a + r_b + r_c)$$
$$= 2\left(|y| + \frac{|-q(y-r)-rx|}{\sqrt{r^2+q^2}} + \frac{|-p(y-r)-rx|}{\sqrt{r^2+p^2}}\right). \quad (7)$$

Let us consider following transformations of equality (6):
$$\sqrt{Q_1} + \sqrt{Q_2} + \sqrt{Q_3} = S,$$
$$\sqrt{Q_1} + \sqrt{Q_2} = S - \sqrt{Q_3},$$
$$(\sqrt{Q_1} + \sqrt{Q_2})^2 = (S - \sqrt{Q_3})^2,$$
$$Q_1 + 2\sqrt{Q_1}\sqrt{Q_2} + Q_2 = S^2 - 2S\sqrt{Q_3} + Q_3,$$
$$2\sqrt{Q_1}\sqrt{Q_2} + 2S\sqrt{Q_3} = S^2 + Q_3 - Q_2 - Q_1,$$
$$(2\sqrt{Q_1}\sqrt{Q_2} + 2S\sqrt{Q_3})^2 = (S^2 + Q_3 - Q_2 - Q_1)^2,$$
$$8S\sqrt{Q_1}\sqrt{Q_2}\sqrt{Q_3} = S^4 - 2S^2 Q_3 - 2S^2 Q_2 - 2S^2 Q_1$$
$$\qquad - 2Q_3 Q_1 - 2Q_3 Q_2 - 2Q_1 Q_2 + Q_3^2 + Q_2^2 + Q_1^2,$$
$$(8S\sqrt{Q_1}\sqrt{Q_2}\sqrt{Q_3})^2 = (S^4 - 2S^2 Q_3 - 2S^2 Q_2 - 2S^2 Q_1$$
$$\qquad - 2Q_3 Q_1 - 2Q_3 Q_2 - 2Q_1 Q_2 + Q_3^2 + Q_2^2 + Q_1^2)^2,$$

after which we obtained an algebraic curve of eight order:
$$64 S^2 Q_1 Q_2 Q_3 = ((S^2 + Q_3 - Q_2 - Q_1)^2 \quad (8)$$
$$\qquad - 4Q_1 Q_2 - 4S^2 Q_3)^2.$$

If we analyze single cases of absolute values in the last equality, in subexpression S given by (7), then we get algebraic equalities of eight order, which fulfill the Erdös-Mordell curve, piece by piece. By application of previous transformation on different ordering of terms $Q_1$, $Q_2$ and $Q_3$, there follows a conclusion that the Erdös-Mordell curve represents the union of parts of algebraic curves of eight order.

The Erdös-Mordell curve gives us a set of telecommunication paths of length $2(r_a + r_b + r_c)$ in a sense of telecommunication problem in [7]. Let us remark that the Erdös-Mordell curve is always in the exterior of the given triangle, and in the special case of equilateral triangle, it additionally contains the center of the triangle as a single isolated point.

### III. NUMERICAL DETERMINING OF $\varepsilon_0$

The open problem given in [6] is proposed as: to prove or disprove that there exists a positive number $\varepsilon_0$ such that the area of $E'$ is larger than $1 + \varepsilon_0$ times the area of the triangle for every triangle in the plane. In [6], there was a conjecture stated: for the finite area of $E'$, the value $\varepsilon_0$ is determined in the case of equilateral triangle $\Delta ABC$. Let us denote:

$P_{\Delta ABC}$ – area of triangle $\Delta ABC$
and

$P_{E'}$ – area of interior parts of the Erdös-Mordell curve.

In this section we describe procedure of numerical determining of constant
$$\varepsilon_0 = \frac{(P_{E'} - P_{\Delta ABC})}{P_{\Delta ABC}} \quad (9)$$
in case of equilateral triangle $\Delta ABC$ with vertices $A=A(0,\sqrt{3})$, $B=B(-1,0)$, $C=C(1,0)$ and area $P_{\Delta ABC} = \sqrt{3}$. Let us form intersection points $P_i = P_i(x_i, y_i)$, $i=1..6$, of the Erdös-Mordell curve with straight lines set through the sides of triangle (Fig. 3):

$$x_1 = -x_4 = -\frac{\sqrt{3}}{\sqrt{15-8\sqrt{3}}} = -1.61966 \ldots$$

$$x_2 = -x_3 = -\frac{1}{2} - \frac{\sqrt{3}}{2\sqrt{15-8\sqrt{3}}} = -1.30983 \ldots$$

$$x_5 = -x_6 = -\frac{1}{2} + \frac{\sqrt{3}}{2\sqrt{15-8\sqrt{3}}} = 0.30983 \ldots$$

$$y_1 = y_4 = 0$$

$$y_2 = y_3 = \frac{\sqrt{3}}{2} - \frac{3}{2\sqrt{15-8\sqrt{3}}} = -0.53664 \ldots$$

$$y_5 = y_6 = \frac{\sqrt{3}}{2} + \frac{3}{2\sqrt{15-8\sqrt{3}}} = 2.26869 \ldots$$

The previous points are obtained by using the symbolic command **solve** of MatLAB.

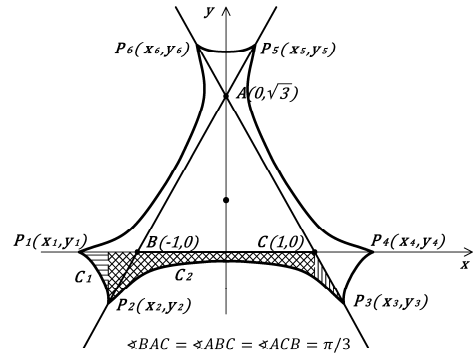

$\varepsilon_0 = 0.8140420779$

Fig. 3. Erdös-Mordell curve (equilateral triangle).

Let us define branch $C_1 = arch\ P_1P_2$ with the appropriate function $y=f_1(x)$ and branch $C_2 = arch\ P_2P_3$ which is defined by the function $y=f_2(x)$. Expressions for $y=f_1(x)$ and $y=f_2(x)$ are determined by using the symbolic command solve of MatLAB, thus solving by $y$ the appropriate algebraic equations of eight order by $x$, piece by piece. Numerical value of constant $\varepsilon_0$ is given by integration:

$$\varepsilon_0 = \sqrt{3}\left(\int_{x_1}^{x_2}(-f_1(x))\,dx + \int_{x_2}^{1}(-f_2(x))\,dx \right.\\ \left. + \int_{1}^{x_3}(-\sqrt{3}x + \sqrt{3} - f_2(x))\,dx\right). \quad (10)$$

Using numerical integration by MatLAB we get:

$$\varepsilon_0 = 0.8140420779\ \ldots \quad (11)$$

Let us note that the previous calculation can be executed in any program that allows solving the polynomial equalities of higher order and numerical integration.

### IV. JAVA APPLET FOR TESTING CONJECTURE AND VISUAL REPRESENTATION OF TRIANGLES FOR GIVEN ε

On the internet at Wolfram Demonstrations Project for The Erdös-Mordell Inequality [9] are available applications that present inner points of Erdös-Mordell curve of the selected triangle. In this paper, the special Java application for visual verification of conjecture from paper [6] will be presented. Let us denote the vertices $A=A(u,v)$, $B=B(-100,0)$ and $C=C(100,0)$ of the triangle $ABC$. Then, for fixed $A$ and arbitrary $M=M(x,y)$ we calculate values:

$$r_a = |y|$$
$$r_b = \frac{|v(x+100) - y(u+100)|}{\sqrt{v^2 + (u+100)^2}}$$
$$r_c = \frac{|v(x-100) - y(u-100)|}{\sqrt{v^2 + (u-100)^2}}$$
$$R_A = \sqrt{(x-u)^2 + (y-v)^2} \quad (12)$$
$$R_B = \sqrt{(x+100)^2 + y^2}$$
$$R_C = \sqrt{(x-100)^2 + y^2}$$
$$g = R_A + R_B + R_C - 2(r_a + r_b + r_c)$$

For values $-600 \le u, v \le 600$ the vertex $A=A(u,v)$ with fixed $B$ and $C$, determines triangle $ABC$ in the visible area of the plane. For every such triangle $ABC$, the equation $g(x,y)=0$ determines the Erdös-Mordell curve and the numerical value of constant $\varepsilon = (P_{E'} - P_{\Delta ABC})/P_{\Delta ABC}$ when Erdös-Mordell curve is closed in domain $-1000 \le x, y \le 1000$. The condition for testing whether the Erdös-Mordell curve is closed is reduced to checking the following inequality $g(x,y) < 0$ for values $x = \pm 1000$, $y = \pm 1000$. Let us note that inequality $g(x,y) \ge 0$ determines all points in domain which are inside of the Erdös-Mordell curve. Since only the integer values with increment of one were used for the points, for each point with positive value of the function $g(x,y)$, value of 1 had been added to area of the Erdös-Mordell curve in square pixels. The integer values were used to improve speed of calculations. As the conjecture was tested on triangles with vertex $A=(u,v)$ where $-600 \le u, v \le 600$, exactly $1200^2 = 1.44 \cdot 10^6$ triangles were tested. Domain of testing points is wider than domain in which triangles are formed, as the Erdös-Mordell curve has at least 81% larger area than the triangle for which it is formed, and in many case much larger. For each triangle, if $-1000 \le x, y \le 1000$, values of $2000^2 = 4 \cdot 10^6$ points were tested for inequality which brings to exactly $1.44 \cdot 10^6 \cdot 4 \cdot 10^6 = 5.76 \cdot 10^{12}$ points processed.

One of the first limitations for application used to prove the conjecture lies in the fact that computers currently lack enough processing power to calculate data, and visually represent them in real time. As a way to circumvent this limitation part of data that would require hours of calculation was done with use of separate application and stored in form that is usable by applications for visualization. To improve the time of calculation, graphical processor unit was used, which allowed massive multithread computation.

The application for visualization of data consists of several areas. The most important area is the one that maps values of ε as equidistant lines. As a visual verification of the conjecture, the color assigned to the value of constant $\varepsilon_0$ appears in the points where the equilateral triangle is formed with two other fixed points, and for $0 \le \varepsilon < \varepsilon_0$ application shows empty screen, as application finds no values that fall in this range. Another area allows visual representation of the Erdös-Mordell curve for user-selected vertexes of triangle. This was achieved by methods of visualization of implicitly expressed surfaces. One of the limitations in visual representation is also the fact that the Erdös-Mordell curve has to be closed. The additional data, such for the area of the Erdös-Mordell curve, are also provided for further analysis in this area. Developed Java application will be available at site http://symbolicalgebra.etf.bg.ac.rs/Java-Applications/.

The conjecture processed in this paper, appertains to experimental mathematics [10]-[12]. Verification of conjecture in one part of visual domain can be expanded and verified also in much wider visual domain.

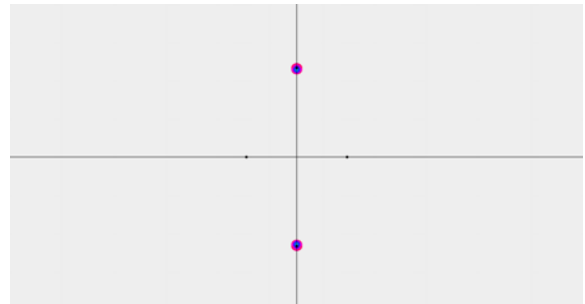

Fig. 4. Color mapping for ε from $\varepsilon_0$ to 0.825.

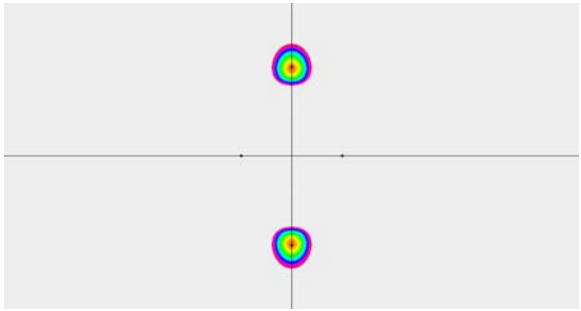

Fig. 5. Color mapping for ε from $\varepsilon_0$ to 1.

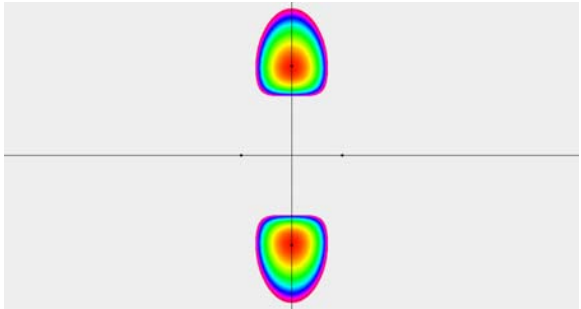

Fig. 6. Color mapping for ε from $\varepsilon_0$ to 2.

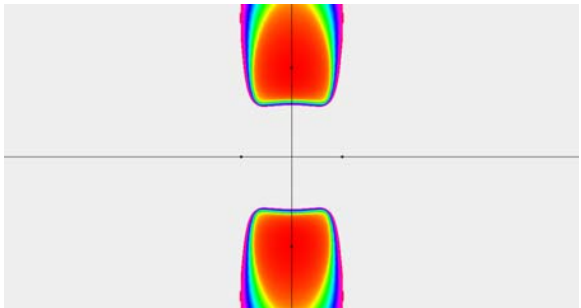

Fig. 7. Color mapping for ε from $\varepsilon_0$ to 18.

## V. CONCLUSION

In the presented paper it is proved that the Erdös-Mordell curve is a union of algebraic curves of eight order and numerical value of constant $\varepsilon_0$ is determined. In addition, Java applet was formed, for visual testing of conjecture that will be the subject of further researches.